\pdfoutput=1

\documentclass[conference]{IEEEtran}
\usepackage[latin9]{inputenc}
\usepackage{amsmath}

\makeatletter
\def\ps@IEEEtitlepagestyle{%
  \def\@oddfoot{\mycopyrightnotice}%
  \def\@evenfoot{}%
}
\def\mycopyrightnotice{%
  {\footnotesize Accepted to IEEE PES GM 2019. \copyright~IEEE 2019 \hfill}
  \gdef\mycopyrightnotice{}
}
\makeatother

\usepackage{amssymb}
\usepackage{amsthm}
\theoremstyle{remark}
\usepackage{booktabs}
\usepackage{color}

\ifCLASSINFOpdf
 \usepackage[pdftex]{graphicx}
 \graphicspath{{figures/}}
\else
\fi
\usepackage{amsmath}
\ifCLASSOPTIONcompsoc
  \usepackage[caption=false,font=normalsize,labelfont=sf,textfont=sf]{subfig}
\else
  \usepackage[caption=false,font=footnotesize]{subfig}
\fi

\usepackage{bbm}

\begin{document}
%
\title{Stochastic Hosting Capacity in \\ LV Distribution Networks}

\author{
\IEEEauthorblockN{Matthew Deakin, \textit{Student Member, IEEE}, Constance Crozier, \textit{Student Member, IEEE},\\ Dimitra Apostolopoulou$^{\dagger}$, \textit{Member, IEEE}, Thomas Morstyn, \textit{Member, IEEE},  \\Malcolm McCulloch, \textit{Senior Member, IEEE}}
\IEEEauthorblockA{Department of Engineering Science, University of Oxford, Oxford, UK\\
$^{\dagger}$Department of Electrical and Electronic Engineering, City University, London, UK}
}


%


\maketitle

\begin{abstract}
Hosting capacity is defined as the level of penetration that a particular technology can connect to a distribution network without causing power quality problems. In this work, we study the impact of solar photovoltaics (PV) on voltage rise. In most cases, the locations and sizes of the PV are not known in advance, so hosting capacity must be considered as a random variable. Most hosting capacity methods study the problem considering a large number of scenarios, many of which provide little additional information. We overcome this problem by studying only cases where voltage constraints are active, with results illustrating a reduction in the number of scenarios required by an order of magnitude. A linear power flow model is utilised for this task, showing excellent performance. The hosting capacity is finally studied as a function of the number of generators connected, demonstrating that assumptions about the penetration level will have a large impact on the conclusions drawn for a given network.
\end{abstract}


%

\section{Introduction}
In this paper we develop a computationally efficient method to calculate the solar photovoltaic (PV) hosting capacity in LV distribution networks. Hosting capacity is typically defined as the total rated power of a given technology that can be connected to a distribution network considering power quality constraints \cite{schwaegerl2005voltage}. If there is uncertainty about the location and number of distributed resources (or of the network conditions, such as load) then the hosting capacity is also uncertain. As such, bounds are determined for hosting capacity \cite{rylander2012stochastic}, which depend on the number and location of generators.

Hosting capacity analysis focuses on a whole range of power quality issues, but one of the most important class of problems are related to voltage quality, particularly voltage rise. In \cite{dubey2015understanding}, the authors run a detailed Monte Carlo simulation to determine hosting capacity, considering PV sizes based on historical PV data. In \cite{dubey2017estimation}, the authors study and illustrate that hosting capacity bounds appear to be roughly distributed according to a Gaussian distribution. In \cite{quijano2017stochastic} the authors also consider the stochasticity of loads and generation to determine an optimal DG penetration. In both \cite{shayani2011photovoltaic} and \cite{jothibasu2016sensitivity}, simplified two-bus models are studied for radial distribution systems, to consider the impact of various network parameters on the hosting capacity of networks. In \cite{ding2016technologies} the authors study the hosting capacity of a number of feeders, considering how the hosting capacity varies if the location of the DG is restricted to specific zones, as well as considering the impact of changes in power factor. In \cite{rylander2016streamlined}, a `streamlined' approach is contrasted with the authors previous `detailed' approach \cite{rylander2012stochastic}. The streamlined method uses network impedance data to study network hosting capacity in a computationally efficient manner. In \cite{navarro2016probabilistic}, the authors study the impact of an increase in low-carbon technologies, including PV, on voltage violations in networks for a range of penetration levels, considering PV sizes distributed according to historical data from the UK.

All of the works surveyed study all scenarios, rather than restricting the studied set to those that have active constraints. This is computationally inefficient, as most scenarios do not improve hosting capacity estimates. Secondly, many methods weight the hosting capacity of 100\% customer penetration equally with penetrations of 1\%. It does not seem to be cost effective (in the long run) for all domestic houses to retrofit to have solar panels, as the costs of ground-based solar PV are approximately half of that of domestic rooftop systems \cite{fu2017us}. As such, weighting in this way is not appropriate.

The contribution of this work is a computationally inexpensive Monte Carlo method for estimating the hosting capacity of distribution networks, choosing only those scenarios for which (voltage) constraints are active. This is based on a linearised power flow, and is compared to a method for which the constraints are not active. Secondly, we study the hosting capacity as an explicit function of the fraction of generators that are connected to the network, rather than concatenating all scenarios into one total hosting capacity. We demonstrate that this leads to considerable differences in the conclusions drawn.

\subsubsection{Notation}
We use $\alpha[k]$ to denote the $k$th element of a vector $\alpha$, $\mathrm{P}(\alpha \,|\, \beta)$ to denote the probability of $\alpha$ given $\beta$, $\mathrm{Re}(\mathbf{z})$ and $\mathrm{Im}(\mathbf{z})$ as the real and imaginary part of a complex quantity $\mathbf{z}$, with bold font representing complex quantities; $\jmath =\sqrt{-1}$ denotes the complex unit, $\mathrm{diag}(\alpha)$ denotes a matrix with the vector $\alpha$ along the diagonal and zeros otherwise, and $\mathbbm{1}^{\alpha} \in \mathbb{R}^{\alpha}$ is a vector of all 1s. Finally, $N_{(\cdot)}$ denotes an integer, while $n_{(\cdot)}\in [0,1]$ denotes a fraction (or percentage).

\section{Voltage-Constrained Hosting Capacity}

\begin{figure}
\centering
\includegraphics[width=0.4\textwidth]{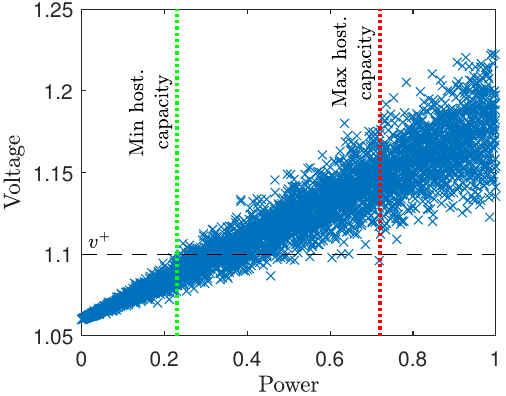}
\caption{The `detailed' hosting capacity method, showing maximum and minimum upper hosting capacities. Note that a large number of scenarios are not close to the constraint and so give little additional information.}\label{f:full}
\end{figure}

For a given network and load, detailed hosting capacity methods consider running a large number of scenarios across a range of PV sizes and locations. The output of the scenarios can then be used to estimate `upper' and `lower' hosting capacities (see Fig. \ref{f:full}). In general, we see that hosting capacity $\Phi $ is a random variable, as we assume that the location of individual PV plants cannot be known a priori. As well as hosting capacity $\Phi$, we also define a power per generator variable $\phi$, given by the relation
\begin{equation}
\phi = \dfrac{\Phi}{N_{\mathrm{gen}}}\,,
\end{equation}
where a feeder has $N_{\mathrm{gen}}$ generators.

In terms of hosting capacity, the optimal location network will usually be at the start of the feeder \cite{rylander2016streamlined}. This is because this is the location in the feeder for which the sensitivity to voltage is smallest. In any case, if the locations of DG are known, then there is no stochastic element and so the hosting capacity is unique.

In this work we assume that all (domestic) loads on the network will wish to connect the same amount of generation. We justify this as
\begin{itemize}
\item the physical availability of space for PV is likely to be relatively similar in small geographical areas (i.e. within an LV network);
\item high levels of socio-economic development tend to be concentrated locally (and thus higher levels of load and abilities to purchase PV), and;
\item government/DSO tariff structures are identical for all customers.
\end{itemize}
Although we enforce this assumption throughout this work, we note that the methods developed could be extended to cases where large fractions of customers will only install PV up to a limit.

We define the $\epsilon $-limited hosting capacity $\Phi_{\epsilon}$ as
\begin{equation}
\mathrm{P}\left(v^{\mathrm{max}} > v^{+} \,\, |\,\, \Phi_{\epsilon}\right) = \epsilon \,, \label{e:hc_definition}
\end{equation}
where $v \in \mathbb{R}^{N_{\mathrm{lds}}}$ is a vector of voltage magnitudes at each of $N_{\mathrm{lds}}$ loads, $v^{\mathrm{max}}$ is the element-wise maximum of $v$, and $v^{+}$ is the voltage upper limit. Setting $\epsilon = 0 \%$ or $\epsilon = 100\%$ correspond to the maximum and minimum hosting capacity.

We assume that the PV generators will be equally likely to connect to any of the loads. Under this assumption, a single scenario consists of a choosing $N_{\mathrm{gen}}$ locations at $N_{\mathrm{lds}}$ load locations. The $i$th scenario is therefore characterised by a set $\Omega_{i}$, which consists of a random subset of load locations. We define the penetration fraction $n_{\mathrm{pen}}$ as
\begin{equation}
n_{\mathrm{pen}} = \dfrac{N_{\mathrm{gen}}}{N_{\mathrm{lds}}}\,.
\end{equation}

\subsection{Approaches to Hosting Capacity}
\subsubsection{The `fixed power' method} The first method we describe as a fixed power hosting capacity approach, and is similar to nominal approaches. A fixed total power $P^{\mathrm{tot}}_{j}$ is chosen and split evenly between $N_{\mathrm{gen}}$ generators, with each load having $P^{\mathrm{gen}}_{j}$ PV generation associated with it as
\begin{equation}
P^{\mathrm{gen}}_{j} = \dfrac{P^{\mathrm{tot}}_{j}}{N_{\mathrm{gen}}}\,.
\end{equation}
The voltage for the $i$th scenario $v_{i}$ is then calculated using
\begin{equation}
v_{i,j} = f_{\mathrm{v}}(P^{\mathrm{gen}}_{j},\Omega_{i})\,,
\end{equation}
where the function $f_{\mathrm{v}}$ describes the mapping from scenario $\Omega_{i}$ with generation $P^{\mathrm{gen}}_{j}$ to voltage magnitudes (i.e. the load flow solution). The probability \eqref{e:hc_definition} is therefore estimated as
\begin{equation}\label{e:n_eps_calc}
\hat{\epsilon}_{j} = \dfrac{\sum_{i=0}^{N_{\mathrm{MC}}-1} \sqcap(v^{\mathrm{max}}_{i,j}> v^{+}) }{N_{\mathrm{MC}}}\,,
\end{equation}
where $\sqcap(\alpha) = 1$ if condition $\alpha$ is true, and is zero otherwise.

To find $\Phi _{\epsilon}$, we use the bisection algorithm, which can be summarised as follows. Two initial guesses for the penetration are chosen ($P_{j=1}^{\mathrm{tot}}$, $P_{j-1=0}^{\mathrm{tot}}$), and $\hat{\epsilon}_{j}$ and  $\hat{\epsilon}_{j-1}$ are calculated according to \eqref{e:n_eps_calc}. If an error function $E$ for two distributions $\hat{\epsilon}_{j}$ and $\hat{\epsilon}_{j-1}$ is greater than some tolerance $\tau$, then a new value $P_{j+1}^{\mathrm{tot}}$ is chosen, according to the sign of $(\hat{\epsilon}_{j} - \epsilon) - (\hat{\epsilon}_{j-1} - \epsilon)$. The process is repeated until convergence is achieved, i.e. until 
\begin{equation}\label{e:tau_defn}
E(\hat{\epsilon}_{j},\hat{\epsilon}_{j-1}) < \tau\,.
\end{equation}
For a more detailed description of the algorithm see, e.g. \cite{gill1981practical4}.

We use the following error metric
\begin{equation}
E(\hat{\epsilon}_{j},\hat{\epsilon}_{j-1}) = \dfrac{|(\hat{\epsilon}_{j} - \epsilon) - (\hat{\epsilon}_{j-1} - \epsilon)|}{1 + |\hat{\epsilon}_{j-1} - \epsilon|}\,,
\end{equation}
which is a `combination' of relative and absolute error \cite{gill1981practical4}.  The algorithm has the advantage that, if ($P_{j=1}^{\mathrm{tot}}$, $P_{j-1=0}^{\mathrm{tot}}$) are either side of a zero, then it is guaranteed to converge.

\subsubsection{The `fixed voltage' method} The second method considers an optimization procedure, which we refer to as a fixed-voltage hosting capacity analysis. In this problem, we only search over cases for which $v^{\mathrm{max}} = v^{+}$. This reduces the number of Monte Carlo runs significantly for the same accuracy. 

To do so, for each scenario $\Omega_{i}$, we solve the optimization problem
\begin{subequations}\label{e:lin_hc_model}
\begin{align}
\hat{P}^{\mathrm{gen}}[i] &= \max \,\, P^{\mathrm{gen}} \\
\mathrm{s.t.} \quad v &\leq v^{+}\mathbbm{1}^{N_{\mathrm{lds}}} \label{e:vplus_constraint}\\
v &= f_{\mathrm{v}}(P^{\mathrm{gen}},\Omega_{i})\,.
\end{align}
\end{subequations}
An estimate for the hosting capacity $\phi _{\epsilon}$ can now be obtained directly from $\hat{P}^{\mathrm{gen}}$, without having to use a method (such as bisection) that was required for the fixed-power approach. 

\begin{figure}
\centering
\includegraphics[width=0.46\textwidth]{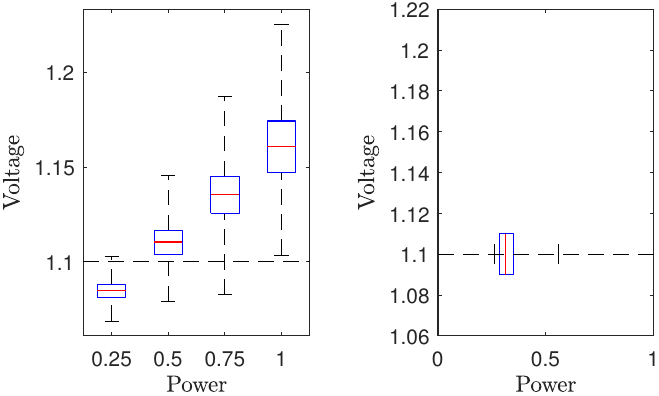}
\caption{The fixed-power hosting capcity method (l), and the fixed-voltage hosting capcity method (r). Enforcing the constraint in the fixed-voltage method reduces the number of scearios required for the same accuracy. (Note: in practise, for the fixed-power method (l), we use a bisection method, rather than the grid-sample method illustrated here, to ensure convergence.)}\label{f:methods}
\end{figure}

A visualization of the two approaches is shown in Fig. \ref{f:methods}. It is clear that the fixed voltage method will require a reduced number of runs for the same accuracy.

\section{Three Phase Linear Power Flow Model}

To reduce the computational complexity, a linear model was developed for the load flow function $f_{\mathrm{v}}$. We utilise a linearization from \cite{bernstein2018load}, which we briefly summarise here. This model has the advantage that it can be linearised around a known load flow solution, increasing the model accuracy, and is valid for three phase unbalanced networks. In general, the method presented would be valid for any three-phase linear model, although the accuracy will vary.

In the case of a three phase, wye-connected distribution network with $N_{\mathrm{bus}}$ 3-phase buses, we can write down the (complex) power flow equations as
\begin{subequations}\label{e:powerflow}
\begin{align}
\mathbf{s} &= \text{diag}(\mathbf{v})\textbf{i}^{*}\,, \\
\mathbf{i} &= \mathbf{Y}\mathbf{v}\,,
\end{align}
\end{subequations}
where $\mathbf{s} \in \mathbb{C}^{3N_{\mathrm{bus}}}$ is the complex power injections at each bus, $\mathbf{i}  \in \mathbb{C}^{3N_{\mathrm{bus}}}$ the corresponding current injections, $\mathbf{Y} \in \mathbb{C}^{3N_{\mathrm{bus}} \times 3N_{\mathrm{bus}}}$ the network admittance matrix and $\mathbf{v} \in \mathbb{C}^{3N_{\mathrm{bus}}}$ the bus complex voltages. 

If we partition the admittance matrix into $\mathbf{Y}_{00}~\in~\mathbb{C}^{3\times 3}, \mathbf{Y}_{L0}~\in~\mathbb{C}^{3(N_{\mathrm{bus}}-1)\times 3}, \mathbf{Y}_{0L}~\in~\mathbb{C}^{3\times 3(N_{\mathrm{bus}}-1)}$, and  $\mathbf{Y}_{LL}~\in~\mathbb{C}^{3(N_{\mathrm{bus}}-1) \times 3(N_{\mathrm{bus}}-1)}$, as
\begin{equation}\label{e:y_partition}
\mathbf{Y} = \begin{bmatrix}
\mathbf{Y}_{00} & \mathbf{Y}_{0L}\\
\mathbf{Y}_{0L} & \mathbf{Y}_{LL}
\end{bmatrix}\,,
\end{equation}
and network voltages and powers as
\begin{equation}\label{e:vs_partition}
\mathbf{v} = \begin{bmatrix}
\mathbf{v}_{0}\\
\mathbf{v}_{L}
\end{bmatrix}\,, \quad \quad
\mathbf{s} = \begin{bmatrix}
\mathbf{s}_{0}\\
\mathbf{s}_{L}
\end{bmatrix}\,.
\end{equation}
where $\mathbf{v}_{0} \in \mathbb{C}^{3}$ is the slack bus voltage, $\mathbf{v}_{L} \in \mathbb{C}^{3(N_{\mathrm{bus}} - 1)}$ the voltages at load buses, $\mathbf{s}_{0} \in \mathbb{C}^{3}$ power injections at the slack bus, and $\mathbf{s}_{L} \in \mathbb{C}^{3(N_{\mathrm{bus}} - 1)}$ the load injections.

From \cite{bernstein2018load}, this implies a model of the form
\begin{equation}\label{e:full_lin_model}
\mathbf{v} = \mathbf{M}
\begin{bmatrix}
p\\
q
\end{bmatrix}
+ \mathbf{a}\,, \quad
\mathbf{M} = 
\begin{bmatrix}
\mathbf{0}_{3\times 3N_{\mathrm{bus}}}\\
\mathbf{M_{\bar{s}}}
\end{bmatrix}\,, \quad
\mathbf{a}= 
\begin{bmatrix}
\mathbf{v}_{0}\\
\mathbf{a_{\bar{s}}}
\end{bmatrix}\,,
\end{equation}
where $\mathbf{s} = p + \jmath q$, and
\begin{subequations}
\begin{align}
\mathbf{M_{\bar{s}}} &= \begin{bmatrix}
\mathbf{Y}_{LL}^{-1}\text{diag}(\bar{\mathbf{v}}_{L})^{-1} & - \jmath \mathbf{Y}_{LL}^{-1}\text{diag}(\bar{\mathbf{v}}_{L})^{-1}
\end{bmatrix}  \\
\mathbf{a_{\bar{s}}} &= - \mathbf{Y}_{LL}^{-1}\mathbf{Y}_{L0}\mathbf{v}_{0}\,,
\end{align}
\end{subequations}
which is linearised at a \textit{known} power flow solution $\bar{\mathbf{v}}$ (for given bus injections $\bar{\mathbf{s}}$).

\subsection{Linear Model Solution}

We now consider the solution of the constrained optimization \eqref{e:lin_hc_model} using the linear model \eqref{e:full_lin_model}. To convert complex voltages $\mathbf{v}$ to voltage magnitudes $v$ we make the assumption that voltage angles do not change significantly from their nominal values. Under this assumption, voltage magnitudes are given by
\begin{equation}\label{e:mag_from_v}
v = \mathrm{Re}\left(\, \mathrm{diag}(\mathbf{v}) \exp (-\jmath \bar{\theta}_{\mathbf{v}}) \,\right) \,,
\end{equation}
where $\bar{\theta}_{\mathbf{v}}$ is vector of the arguments of the complex voltages at the known load flow solution $\mathbf{\bar{v}}$. As such, we write done
\begin{equation}
v = F\begin{bmatrix}
p\\
q
\end{bmatrix} + g\,,
\end{equation}
where
\begin{align}
F &= \begin{bmatrix}
\mathrm{diag}(\cos(\bar{\theta}_{\mathbf{v}})) & \mathrm{diag}(\sin(\bar{\theta}_{\mathbf{v}}))
\end{bmatrix}\begin{bmatrix}
\mathrm{Re}(\mathbf{M})\\
\mathrm{Im}(\mathbf{M})
\end{bmatrix} \,,\\
g &= \begin{bmatrix}
\mathrm{diag}(\cos(\bar{\theta}_{\mathbf{v}})) & \mathrm{diag}(\sin(\bar{\theta}_{\mathbf{v}}))
\end{bmatrix}\begin{bmatrix}
\mathrm{Re}(\mathbf{a})\\
\mathrm{Im}(\mathbf{a})
\end{bmatrix}\,.
\end{align}

The load vector $\mathbf{s}$ consists of the fixed (linearised) component, plus the generation component, $\mathbf{s} = \bar{\mathbf{s}} + P^{\mathrm{gen}}\Lambda$, where the vector $\Lambda$ is given as
\begin{align*}
\Lambda [k] = 
\begin{cases}
1 \quad \mathrm{if }\; k \in \Omega_{i} \\
0 \quad \mathrm{otherwise}\,.
\end{cases}
\end{align*}
Noting that we have just one unknown, $P^{\mathrm{gen}}$, we can write down the the voltage at each node as 
\begin{equation}\label{e:final_lin_model}
v = F\Lambda P^{\mathrm{gen}} + \bar{v}\,,
\end{equation}
Therefore, the solution to \eqref{e:lin_hc_model} is found (by substituting \eqref{e:final_lin_model} into \eqref{e:vplus_constraint}) as
\begin{equation}
\hat{P}^{\mathrm{gen}} = \min (\lfloor \mathrm{diag}(F\Lambda)^{-1}(\mathbbm{1}^{N_\mathrm{lds}}v^{+} - \bar{v}) \rfloor)\,,
\end{equation}
where we have used the notation 
\begin{equation}
\lfloor \alpha[k] \rfloor = \begin{cases}
\alpha[k] \quad \mathrm{if} & \alpha[k] > 0 \\
\infty & \mathrm{otherwise}\,.
\end{cases}
\end{equation}
This optimization is therefore very fast to run as it only consists of arithmetic operations and sorting.

\textit{Remark: Non-uniform generator selection.} The model we have used here uses the approximation that all houses will wish to install as much PV as possible. In the case that there is information to suggest that this will not be the case (i.e. that a significant fraction of houses will want less than the maximum PV), then additional variables can be added. This would result in a linear program, for which there are efficient solvers available.

\section{Case Studies}

In this work we study four low voltage networks from \cite{enwl2015low} and the European low voltage network \cite{ieee2017distribution} (see Table~\ref{t:network_details}). These networks are notable in that they do not have any voltage control equipment (capacitors or voltage regulators). The linear model is therefore particularly well suited to modelling these types of networks. Note that model 1 and model 2 are identical, save for a change in substation voltage. Each load is assigned a demand of 0.3 kW at 0.95 power factor lagging.


\begin{table}
\caption{Feeder details (NX.Y as network X, feeder Y)}\label{t:network_details}
\centering
\begin{tabular}{llllll}
\toprule
 & EU LV & N1.1 & N2.1 & N3.1 & N4.1 \\
\midrule
No. loads $N_{\mathrm{lds}}$  & 55 & 55 & 175 & 94 & 24\\
No. buses $N_{\mathrm{bus}}$ & 907 & 907 & 2287 & 1304 & 375\\
Substation voltage $v^{\mathrm{Sub}}$, pu & 1.05 & 1.00 & 1.00 & 1.00 & 1.00\\
\bottomrule
\end{tabular}
\end{table}


\subsection{Validation and comparison of proposed approaches}

\subsubsection{Validation of linear model} 
\begin{figure}
\centering
\includegraphics[width=0.4\textwidth]{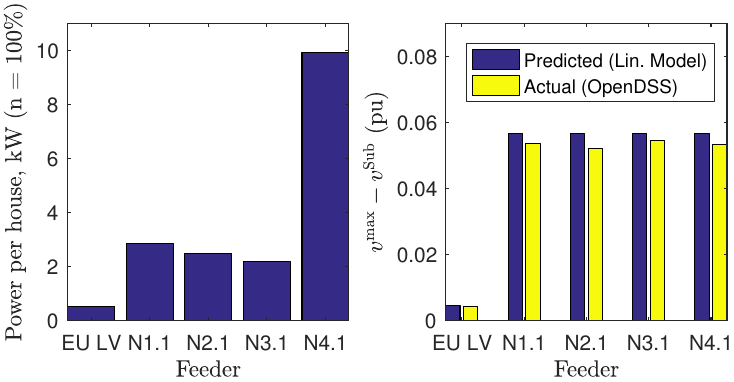} 
\caption{Power that can be exported per house versus feeder, at 100$\%$ penetration (l), and the predicted and actual voltage rise on the network ($v^{\mathrm{Sub}}$ is the substation voltage) (r).}\label{f:lin_program_vldt}
\end{figure}

To first validate the accuracy of the linear model for this model, the difference between the actual and predicted maximum voltage rise was calculated, for the case of 100\% penetration (for which the hosting capacity is deterministic). It was found that the linear model behaves well compared to a full load flow solution (see Fig. \ref{f:lin_program_vldt}), with the non-linear load flow calculated in OpenDSS \cite{opendss2017}).

\subsubsection{Monte Carlo simulation validation} The number of runs required for the monte-carlo method $N_{\mathrm{MC}}$ to achieve good accuracy was tested considering a hosting capacity $\Phi_{\epsilon}$ with $\epsilon = 5\%$ (that is, if $P^{\mathrm{tot}} = \Phi _{5\%}$ then there is a 5\% chance of there being a steady state overvoltage). The calculated hosting capacity for two Monte Carlo runs for the fixed-voltage method is given in Table \ref{tab:montecarlo_validation}. We see that using $N_{\mathrm{MC}}=1000$ gives an accuracy better than 3\% in the networks studied.

\begin{table}
	\begin{center}
	\caption{Estimated hosting capacity and error for two monte carlo runs ($N_{\mathrm{MC}}$ = 1000, $n_{\mathrm{pen}}=50\%$)}
	\label{tab:montecarlo_validation}
		\begin{tabular}{rlll}
		\toprule
			{Feeder} & {$\Phi_{5\%}$ (run A), kW} & {$\Phi_{5\%}$ (run B), kW} & {Rel. error, \%}\\
			\midrule
					{EU LV} & $15.3$ & $14.9$ & $2.27$\\
					{N1.1} & $92.4$ & $89.9$ & $2.68$\\
					{N2.1} & $230.1$ & $229.7$ & $0.15$\\
					{N3.1} & $115.5$ & $116.3$ & $0.71$\\
					{N4.1} & $120.5$ & $122.6$ & $1.70$\\
		\bottomrule
		\end{tabular}
	\end{center}
\end{table}

\subsubsection{Comparison of fixed power and fixed voltage methods} A comparison of the computational efficiency of the fixed-power and fixed-voltage methods are given in Table \ref{tab:montecarlo_validation}. In the fixed power method, we use a tolerance $\tau = 1\%$ in \eqref{e:tau_defn}, initialised with $P_{0}^{\mathrm{tot}}=0$, and $P_{1}^{\mathrm{tot}}$ as the unique hosting capacity at $n_{\mathrm{pen}}=100\%$.

The number of iterations required to reach convergence of the bisection method is clearly a key driver in the increase in time required to run the fixed-power method compared to the fixed-voltage method. The difference in accuracy is of the same order of magnitude as the difference caused by running multiple Monte Carlo runs. Given the computational speed, we conclude that the fixed voltage method has a clear advantage over the fixed power approach.

\begin{table}
	\begin{center}
	\caption{Comparison of timings and estimated hosting capacities for the fixed power and fixed voltage methods}
	\label{tab:montecarlo_comparison}
		\begin{tabular}{rlllll}
\toprule
& \multicolumn{3}{c}{Fixed power} & \multicolumn{2}{c}{Fixed voltage} \\
\cmidrule(l{0.3em}r{0.6em}){2-4}\cmidrule(l{0.3em}r{0.6em}){5-6}
			{Feeder} & {Iterations} & {Time, s} & {$\Phi _{5\%}$, kW} & {Time, s} & {$\Phi _{5\%}$, kW}\\
			\midrule
					{EU LV} & $8$ & $5.30$ & $15.0$ & $0.80$ & $15.3$\\
					{N1.1} & $11$ & $6.95$ & $91.4$ & $0.77$ & $92.4$\\
					{N2.1} & $23$ & $102.83$ & $228.4$ & $5.57$ & $230.1$\\
					{N3.1} & $10$ & $15.49$ & $115.2$ & $1.82$ & $115.5$\\
					{N4.1} & $10$ & $1.20$ & $120.2$ & $0.17$ & $120.5$\\
		\bottomrule
		\end{tabular}
	\end{center}
\end{table}

\subsection{Hosting Capacity for DG Policy Decision}

\begin{figure}
\centering
\includegraphics[width=0.45\textwidth]{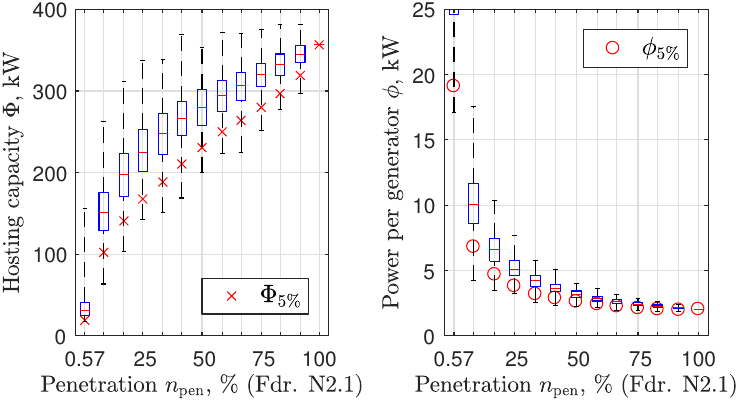}
\caption{Boxplot of the hosting capacity $\Phi$ (l) and power per generator $\phi$ (r) as a function of the penetration $N_{\mathrm{gen}}$ for feeder N2.1.}\label{f:variable_N_3}
\end{figure}

\begin{figure}
\centering
\includegraphics[width=0.45\textwidth]{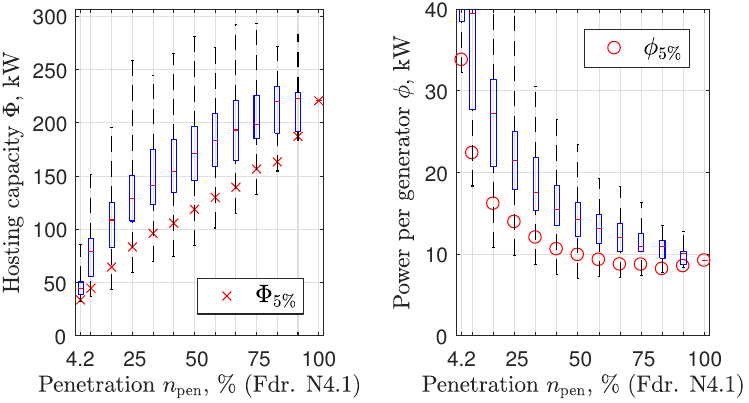}
\caption{Boxplot of the hosting capacity $\Phi$ (l) and power per generator $\phi$ (r) as a function of the penetration $N_{\mathrm{gen}}$ for feeder N4.1.}\label{f:variable_N_5}
\end{figure}

We now study hosting capacity as a function of the customer penetration level. The probability density function (PDF) of the power per house and total feeder power is shown in Fig. \ref{f:variable_N_3} for feeder N2.1 and in Fig. \ref{f:variable_N_5} for feeder N4.1. The 5\% hosting capacity $\Phi _{5\%}$ is plotted alongside boxplots of the estimate of the hosting capacity, illustrating the minimum and maximum, interquartile range, and median of $\Phi$; the same is also plotted for $\phi$. 

It is clear that the hosting capacity is strongly affected by the fraction of loads connected. In both feeders, the 5\% hosting capacity $\Phi _{5\%}$ appears to loosely follow an S-shaped curve: it increases rapidly first, then increasing at a linear rate, before increasing more rapidly towards 100\% penetration. 

The extrema (the minimum and maximum hosting capacity) appear to have some noise. That is, it seems that a larger number of Monte Carlo runs would be required to estimate these accurately. The minimum hosting capacity largely follows the S-shaped curve of the 5\% hosting capacity, while maximum hosting capacity follows an inverted-U shaped curve, increasing rapidly at first before levelling off, finally decreasing slightly at 100\% penetration.

The amount of power that can be connected per generator, $\phi$, also varies significantly. In both networks we see that that median power per generator drops monotonically. On the other hand, the minimum and 5\% power per generator reach a minimum close to 80\% for feeder N4.1 (see Fig. \ref{f:variable_N_5}, (r)). This is presumed to be due to the possibility of greater unbalance with smaller numbers of generators.

We note that the power per generator at a penetration of 25\% is 50-100\% greater than at 100\% penetration in both networks. If DSOs assume that only small numbers of households will connect (for example, if a community develops a larger community-owned PV scheme), then they can allow more PV to connect. This increases the utilisation of the network. Finally, although they have not been studied in this work, we note that other power quality constraints will also impact on the results (for example, thermal constraints).

\section{Conclusions}
Hosting capacity of domestic PV will always have a stochastic element to it, and finding computationally efficient, transparent methods of calculating PV hosting capacity is paramount given the large number of LV feeders. The fixed-voltage hosting capacity method will equip DSOs with the ability to rapidly study hosting capacity in a computationally efficient manner.

Given recent policy shifts towards larger-scale PV, in future it may be that relatively fewer PV systems connect, but the PV systems that connect are larger. We have demonstrated that the total hosting capacity is a strong function of the number of PV generators that are connected. Depending on the expected uptake of PV, the methods present will help with decision making, evaluating the hosting capacity in a fair, transparent way.

\section*{Acknowledgement}

The authors wish to thank the Oxford Martin Programme on Integrating Renewable Energy, the John Aird Scholarship, and the Clarendon Scholarship for their support.



%
\bibliographystyle{IEEEtran}
\bibliography{pesgm19_bib}{}

\end{document}